\documentclass[a4paper,12pt]{article}

\usepackage{amsmath}
\usepackage{amsthm}
\usepackage{amssymb}
\newtheorem{thm}{Theorem}
\newtheorem{prop}[thm]{Proposition}

\newtheorem{cor}[thm]{Corollary}
\newtheorem{rem}{Remark}

\theoremstyle{definition}
\newtheorem{defn}[thm]{Definition}

  \newtheorem{df-thm}[thm]{Definition-Theorem}
  \newtheorem{df-prop}[thm]{Definition-Proposition}
  \newtheorem{df-lem}[thm]{Definition-Lemma}
  \newtheorem{df-exa}[thm]{Definition-Example}
  
  \newtheorem{exa}[rem]{Example}
\newcommand{\adag}{a^{\dagger}}
\newcommand{\tr}{{\rm tr}}
\title{
Noncommutative Spectral Decomposition with Quasideterminant}
\author{Tatsuo Suzuki\thanks{
  {\it E-mail address}: suzukita{\char'100}aoni.waseda.jp}\\
\\
  Department of Mathematical Sciences, 
  Waseda University, \\ 
  Tokyo 169-8555, 
  Japan}
\date{}
\begin{document}
\maketitle

\begin{abstract}
We develop a noncommutative analogue of the spectral decomposition 
with the quasideterminant defined by I. Gelfand and V. Retakh. 
In this theory, 
by introducing a noncommutative Lagrange interpolating polynomial and 
combining a noncommutative Cayley-Hamilton's theorem and 
an identity given by a Vandermonde-like quasideterminant, we can 
systematically calculate a function of a matrix even if it has 
noncommutative entries. 
As examples, the noncommutative spectral 
decomposition and the exponential matrices of a quaternionic matrix and of 
a matrix with entries being harmonic oscillators are given. 
\end{abstract}

\newpage

\section{Introduction}

The theory of spectral decomposition of a square matrix over 
a commutative field is well-known in linear algebra and is used for 
calculation of a function of the matrix, especially the exponential matrix. 
However, for a matrix with noncommutative entries, 
the determinant or the characteristic polynomial are not defined 
because of the ordering problem. Therefore, ``eigenvalues" used in 
the spectral decomposition are undefined and we have no systematic method 
for calculation of function of a matrix with noncommutative entries 
until now. 

Under these circumstances, we studied the exponential of a matrix 
with entries being harmonic osillators for a model in quantum optics and 
developed ``the quantum diagonalization method" for a special type of 
matrices derived from the representation theory \cite{FHKSW}. Moreover, 
we had a chance to encounter with the quasideterminant defined by 
I. Gelfand and V. Retakh. By using the quasideterminant, 
``noncommutative determinants" such as quaternionic determinants \cite{As}, 
superdeterminant, quantum determinant, Capelli determainant, etc. are 
expressed in the unified form \cite{GR1}. In the theory of the noncommutative 
integrable system, quasideterminants are very useful 
to express the solution of the noncomutative integrable equations 
\cite{EGR}, \cite{H}, \cite{GN}, \cite{GNO}. Furthermore, 
various noncommutative analogue of theories using determinants are developed, 
for example, noncommutative analogue of Cramer's formula, the Vandermonde 
determinant, symmetric functions, Pl\"ucker coordinates, and so on. 
(see, \cite{GR2}, 
\cite{GGRW}, \cite{GKLLRT} and references within). 

In particular, in \cite{GKLLRT}, they investigated a noncommutative 
Cayley-Hamilton's theorem. In their theory, a different 
characteristic polynomial for each row was introduced and 
the trace or determinant were of the form of diagonal matrices. 
Moreover, we knew through the study of the quantum diagonalization method 
that eigenvalues should be generalized as ``eigen-diagonalmatrics" 
due to the noncommutativity of entries of the matrix. 
That is why we find that a noncommutative Cayley-Hamilton's theorem 
in \cite{GKLLRT} is suitable to a noncommutative analogue of 
the spectral decomposition. 

In this paper, we define a noncommutative analogue of the 
Lagrange interpolating polynomial and develop a noncommutative analogue of 
the spectral decomposition by using 
the noncommutative Cayley-Hamilton's theorem with the quasideterminant. 
An identity given by a Vandermonde-like quasideterminant 
plays an essential role. 
As examples, we explicitly calculate the noncommutative spectral 
decomposition and the exponential matrices of a quaternionic matrix and of 
a matrix with entries being harmonic oscillators. 

The contents of this paper are as follows. 
In section 2, we give a brief review of the spectral decomposition in 
linear algebra. In section 3, we introduce the quasideterminant 
defined by I. Gelfand and V. Retakh and describe some important properties 
used in our theory. In section 4, we review the noncommutative 
Cayley-Hamilton's theorem in \cite{GKLLRT} shortly. In section 5, 
we develop a noncommutative analogue of the spectral decomposition 
with the quasideterminant. In section 6, we apply our method to 
a quaternionic matrix and a matrix with entries being 
harmonic oscillators. 
Section 7 is devoted to disscution. 

\newpage

\section{Brief Review of the Spectral Decomposition}

Firstly, we give a brief review of the spectral decomposition in 
linear algebra. 

Let $A$ be a $n \times n$-matrix with commutative entries. 
For simplicity, we suppose that all the eigenvalues 
$\lambda_1, \cdots, \lambda_n$ of $A$ are distinct. 
For $j=1,\cdots,n$, we set 
$$P_j=\prod_{1 \le i \le n, \ i \ne j }
\frac{(A-\lambda_i I)}{(\lambda_j-\lambda_i)}. $$
The polynomial of right hand side is called the Lagrange interpolating 
polynomial. Then we have the spectral decomposition of $A$; 
$$A=\lambda_1 P_1+\cdots+ \lambda_n P_n. $$
Moreover, if {\underline{the Cayley-Hamilton's theorem holds}}，then 
$P_1, \ \cdots, \ P_n \ $ are projection matrices，i.e.
$$P_i^2=P_i, \quad P_iP_j=O \ (i \ne j), \quad 
P_1+\cdots+P_n=I. $$
Therefore, we can calculate $\ \exp A \ $ explicitly; 
$$\exp A=e^{\lambda_1} P_1+\cdots+ e^{\lambda_n} P_n. $$

\begin{rem}
Lagrange interpolating polynomials 
$\displaystyle \ f_j(z)=\prod_{1 \le i \le n, \ i \ne j }
\frac{(z-x_i)}{(x_j-x_i)} \quad (j=1,\cdots,n)\ $ satisfy the following 
relations; \\
$(1) \quad x_1^j f_1(z)+x_2^j f_2(z)+\cdots+x_n^j f_n(z)=z^j 
\quad (j=0,1,\cdots,n-1)$\\
$(2) \quad f_i(x_j)=\delta_{ij}$\\
We note that from (1), if $x_1,\cdots,x_n$ are all distinct, then we have 
$$\begin{pmatrix}
f_1(z) \\
\vdots \\
f_n(z)
\end{pmatrix}
=\begin{pmatrix}
x_1^{n-1} & \cdots & x_n^{n-1} \\
& \cdots & \\
1 & \cdots & 1
\end{pmatrix}^{-1}
\begin{pmatrix}
z^{n-1} \\
\vdots \\
1
\end{pmatrix}. 
$$
\end{rem}

\section{Quasideterminant}
\label{sec:pre}

In this section, we introduce the quasideterminant 
defined by I. Gelfand and V. Retakh and describe some important properties 
used in our theory. 

\subsection{Definition}

Let $R$ be a (not necessary commutative) associative algebra. 
For a position $(i,j)$ in a square matrix 
$A=(a_{rs})_{1 \le r,s \le n} \in M(n,R)$, let 
$A^{ij}$ denote the $(n-1)\times(n-1)$-matrix obtained from $A$ by deleting 
the $i$-th row and the $j$-th column. Let also 
$\mbox{r}_i{}^j=(a_{i1},\cdots,\hat{a}_{ij},\cdots,a_{in})$ and 
$\mbox{c}_j{}^i=(a_{1j},\cdots,\hat{a}_{ij},\cdots,a_{nj})^T$. 


\begin{defn}
We assume that $A^{ij}$ is invertible over $R$. 
The $(i,j)$-quasideterminant of $A$ is defined by
\begin{equation}
|A|_{ij}=a_{ij}-\mbox{r}_i{}^j \cdot (A^{ij})^{-1} \cdot \mbox{c}_j{}^i. 
\end{equation}
\end{defn}

\noindent
\begin{exa}
For $A=\left( \begin{array}{cc}
a_{11} & a_{12} \\
a_{21} & a_{22}
\end{array} \right)$, 
$$ |A|_{11}=a_{11}-a_{12}a_{22}^{-1}a_{21},
\quad |A|_{12}=a_{12}-a_{11}a_{21}^{-1}a_{22}, $$
$$ |A|_{21}=a_{21}-a_{22}a_{12}^{-1}a_{11},
\quad |A|_{22}=a_{22}-a_{21}a_{11}^{-1}a_{12}.$$
\end{exa}

\noindent
It is sometimes convinient to adopt the following more explicit notation
$$ |A|_{11}=
\left| 
\begin{array}{cc}
\cline{1-1}
\multicolumn{1}{|c|}{a_{11}} & a_{12} \\ \cline{1-1}
a_{21} & a_{22} 
\end{array}
\right|
=a_{11}-a_{12}a_{22}^{-1}a_{21}. $$

\begin{rem}
If the elements $a_{ij}$ of the matrix $A$ commute, then
$$|A|_{ij}=(-1)^{i+j}\frac{\det A}{\det A^{ij}}.$$
\end{rem}

\subsection{Inverse Matrix and Quasideteminant}

\begin{prop}
If all $|A|_{ij}^{-1}$ exist, $A^{-1}$ is given by 
$$A^{-1}=(|A|_{ji}^{-1})_{1 \le i,j \le n}. $$
\end{prop}

\begin{exa}
For a quaternionic matrix $A=\left( \begin{array}{cc}
1 & i \\
j & k 
\end{array} \right)$, \\
\begin{eqnarray*}
|A|_{11}^{-1}&=&(1-i\cdot k^{-1} j)^{-1}=(1+ikj)^{-1}=\frac{1}{2}\\
|A|_{21}^{-1}&=&(j-k\cdot i^{-1} 1)^{-1}=(j+ki)^{-1}=(2j)^{-1}
 =-\frac{j}{2}\\
|A|_{12}^{-1}&=&(i-1\cdot j^{-1} k)^{-1}=(i+jk)^{-1}=(2i)^{-1}
 =-\frac{i}{2}\\
|A|_{22}^{-1}&=&(k-j\cdot 1^{-1} i)^{-1}=(k-ji)^{-1}=(2k)^{-1}
 =-\frac{k}{2}. 
\end{eqnarray*}
Therefore 
$$A^{-1}=\frac{1}{2} \left( \begin{array}{cc}
1 & -j \\
-i & -k 
\end{array} \right) .$$
\end{exa}

\begin{exa}\ We can calculate quasideterminants inductively;
\begin{eqnarray*}
&& \hspace{-1cm}
\left| \begin{array}{ccc}
\cline{1-1}
\multicolumn{1}{|c|}{a_{11}}& a_{12} & a_{13} \\ \cline{1-1}
a_{21} & a_{22} & a_{23} \\
a_{31} & a_{32} & a_{33} \\
\end{array} \right| 
=a_{11}
-(a_{12} \ a_{13})
\left( \begin{array}{cc}
|A^{11}|_{22}^{-1} & |A^{11}|_{32}^{-1} \\
|A^{11}|_{23}^{-1} & |A^{11}|_{33}^{-1} 
\end{array} \right)
\left( \begin{array}{c}
a_{21} \\
a_{31}
\end{array} \right) \\
&& \\ 
&=&a_{11}
-a_{12}(a_{22}-a_{23}a_{33}^{-1}a_{32})^{-1}a_{21}
-a_{12}(a_{32}-a_{33}a_{23}^{-1}a_{22})^{-1}a_{31}\\ && \quad \ \ 
-a_{13}(a_{23}-a_{22}a_{32}^{-1}a_{33})^{-1}a_{21}
-a_{13}(a_{33}-a_{32}a_{22}^{-1}a_{23})^{-1}a_{31}.
\end{eqnarray*}
\end{exa}

\subsection{Homological Relations}

For $A=(a_{ij}) \in M(n,R)$, $n^2$ quasideterminants are defined. 
They are related by the so-called homological relations. For example, 
$$
\left| \begin{array}{cc}
a_{11} & a_{12} \\ \cline{1-1}
\multicolumn{1}{|c|}{a_{21}} & a_{22} \\ \cline{1-1}
\end{array} \right| 
=-a_{22} \ a_{12}^{-1} 
\left| \begin{array}{cc}
\cline{1-1}
\multicolumn{1}{|c|}{a_{11}} & a_{12} \\ \cline{1-1}
a_{21} & a_{22}
\end{array} \right|.
$$
In general, we have important identities as follows; 
\begin{prop}
\begin{enumerate}
\item Row homological relations: 
$$ -|A|_{ij} \cdot |A^{il}|_{sj}^{-1}
=|A|_{il} \cdot |A^{ij}|_{sl}^{-1}, \qquad s \ne i$$
\item Column homological relations: 
$$ -|A^{kj}|_{it}^{-1} \cdot |A|_{ij}
=|A^{ij}|_{kt}^{-1} \cdot |A|_{kj}, \qquad t \ne j$$
\end{enumerate}
\end{prop}

\subsection{The Multiplication of Rows and Columns}

Let $B$ be the matrix obtained from the matrix $A$ by multiplying 
the $i$-th row by $\lambda \in R$，then 
$$|B|_{kj}=\left\{ 
  \begin{array}{ccc}
  \lambda |A|_{ij} & \mbox{if} & k=i \\
  |A|_{kj}  & \mbox{if} & k \ne i 
  \end{array} \right.
$$
Let $C$ be the matrix obtained from the matrix $A$ by multiplying 
the $j$-th column by $\mu \in R$，then 
\begin{equation}
|C|_{il}=\left\{ 
  \begin{array}{ccc}
  |A|_{ij} \mu & \mbox{if} & l=j \\
  |A|_{il}  & \mbox{if} & l \ne j 
  \end{array} \right.
\label{eqn:mu}
\end{equation}

\begin{exa}
$$\left| 
\begin{array}{cc}
\cline{2-2}
a_{11} & \multicolumn{1}{|c|}{a_{12}\mu} \\ \cline{2-2}
a_{21} & a_{22}\mu
\end{array}
\right|
=a_{12}\mu -a_{11}a_{21}^{-1}a_{22}\mu=|A|_{12}\mu
$$
$$\left| 
\begin{array}{cc}
\cline{2-2}
a_{11}\mu & \multicolumn{1}{|c|}{a_{12}} \\ \cline{2-2}
a_{21}\mu & a_{22}
\end{array}
\right|
=a_{12}-a_{11}\mu (a_{21}\mu)^{-1}a_{22}=|A|_{12}
$$
\end{exa}

\subsection{Sylvester's Identity}

Let $A=(a_{ij}) \in M(n,R)$ be a matrix 
and $A_0=(a_{ij}), i,j=1,\cdots,k$, a submatrix of $A$ that is invertible 
over $R$. For $p,q=k+1,\cdots,n$, set 
\begin{equation}
c_{pq}=
\left| \begin{array}{cccc}
& & & a_{1q} \\
& A_0 & & \vdots \\
& & & a_{kq} \\
\cline{4-4}
a_{p1} & \cdots & a_{pk} & \multicolumn{1}{|c|}{a_{pq}} \\
\cline{4-4}
\end{array} \right| .
\label{eqn:pivot}
\end{equation}
These quasideterminants are defined because matrix $A_0$ is invertible. 

Consider the $(n-k) \times (n-k)$ matrix
$$C=(c_{pq}), \quad p,q=k+1,\cdots,n.$$
The submatrix $A_0$ is called the {\it{pivot}} for the matrix $C$. 
\begin{thm}(Sylvester's identity) \ 
For $i,j=k+1,\cdots,n$, 
$$|A|_{ij}=|C|_{ij}.$$
\label{thm:Syl}
\end{thm}

\begin{exa}
$$
\left|
\begin{array}{ccc}
1 & a_{12} & a_{13} \\
0 & a_{22} & a_{23} \\
0 & a_{32} & a_{33} \\
\end{array}
\right|_{ij} \quad (i,j=2,3)
$$
Applying the theorem \ref{thm:Syl} with 
the $(1,1)$-entry $1$ as a pivot, we put 
$$c_{pq}=\left|
\begin{array}{ccc}
1 & a_{1q} \\
0 & a_{pq} \\
\end{array}
\right|_{pq}=a_{pq} \quad (p,q=2,3)
$$
and 
$$
\left|
\begin{array}{ccc}
1 & a_{12} & a_{13} \\
0 & a_{22} & a_{23} \\
0 & a_{32} & a_{33} \\
\end{array}
\right|_{ij}
=|C|_{ij}=
\left|
\begin{array}{ccc}
a_{22} & a_{23} \\
a_{32} & a_{33} \\
\end{array}
\right|_{ij} \quad (i,j=2,3).
$$
\end{exa}

\section{Noncommutative Version of the Characteristic Polynomial and 
the Cayley-Hamilton's Theorem}

In this section, we review the noncommutative Cayley-Hamilton's 
theorem in \cite{GKLLRT} 
shortly. We use notations $\Phi_i(\lambda), \ C_{(i)j}$ instead of 
$Q_i(t), \ L_{j}^{(i)}(A)$ in it. 
For $\ \displaystyle A=
\begin{pmatrix}
a_{11} & a_{12} \\
a_{21} & a_{22} 
\end{pmatrix}$
, we denote $\Phi_1(\lambda), \Phi_2(\lambda)$ as two polynomials 
given by 
\begin{eqnarray*}
\Phi_1(\lambda)&=&
\lambda^2-(a_{11}+a_{12}a_{22}a_{12}^{-1})\lambda+
(a_{12}a_{22}a_{12}^{-1}a_{11}-a_{12}a_{21})\\
&\equiv&
\lambda^2-\tr_1(A)\lambda+\det{}_1 (A),
\end{eqnarray*}
\begin{eqnarray*}
\Phi_2(\lambda)&=&
\lambda^2-(a_{22}+a_{21}a_{11}a_{21}^{-1})\lambda+
(a_{21}a_{11}a_{21}^{-1}a_{22}-a_{21}a_{12})\\
&\equiv&
\lambda^2-\tr_2(A)\lambda+\det{}_2 (A).
\end{eqnarray*}
Then we can check the noncommutative Cayley-Hamilton's theorem for 
the generic matrix of order $2$; 
$$A^2-
\begin{pmatrix}
\tr_1(A) & 0 \\
0 & \tr_2(A)
\end{pmatrix}A
+\begin{pmatrix}
\det{}_1(A) & 0 \\
0 & \det{}_2(A)
\end{pmatrix}=O.
$$
The general result is as follows. We also give a simple proof. 
\begin{thm}\cite{GKLLRT}
For $A=(a_{ij})\ \in M(n,R)$, we define a ``noncommutative characteristic 
polynomial for the $i$-th row" as follows; \\
\begin{eqnarray}
\Phi_i(\lambda)&=&
\left| 
\begin{array}{ccccc}
\cline{5-5}
a_{i1}^{(n)} & a_{i2}^{(n)} & \cdots & a_{in}^{(n)} & 
\multicolumn{1}{|c|}{\lambda^n} \\ 
\cline{5-5}
a_{i1}^{(n-1)} & a_{i2}^{(n-1)} & \cdots & a_{in}^{(n-1)} & \lambda^{n-1}\\
\vdots & & & \vdots & \vdots \\
a_{i1}^{(1)} & a_{i2}^{(1)} & \cdots & a_{in}^{(1)} & \lambda \\ 
a_{i1}^{(0)} & a_{i2}^{(0)} & \cdots & a_{in}^{(0)} & 1 \\
\end{array}
\right| \label{eqn:NCCP}\\
&\equiv &\lambda^n-\sum_{k=1}^n C_{(i)k}\lambda^{n-k}, \nonumber
\end{eqnarray}
where $A^k=(a_{ij}^{(k)})$. Then we have a noncommutative version of 
the Cayley-Hamilton theorem 
\begin{equation}
A^n-\sum_{k=1}^n
\begin{pmatrix}
C_{(1)k} & & & \\
& C_{(2)k} & & \\
& & \ddots & \\
& & & C_{(n)k}
\end{pmatrix}
A^{n-k}=O. \label{eqn:NCCH}
\end{equation}
\end{thm}
\begin{proof}
For unknown $C_{(i)k} \ (i,k=1,\cdots,n)$, consider the equation 
(\ref{eqn:NCCH}). Then the $(i,j)$-entry of (\ref{eqn:NCCH}) is 
\begin{equation}
a_{ij}^{(n)}-\sum_{k=1}^n C_{(i)k}a_{ij}^{(n-k)}=0 \quad 
(i,j=1,\cdots,n), 
\label{eqn:coef1}
\end{equation}
namely 
\begin{equation}
(C_{(i)1},\cdots,C_{(i)n})
\begin{pmatrix}
a_{i1}^{(n-1)} & \cdots & a_{in}^{(n-1)} \\
\vdots & & \vdots \\
a_{i1}^{(0)} & \cdots & a_{in}^{(0)}
\end{pmatrix}
=(a_{i1}^{(n)},\cdots,a_{in}^{(n)}). 
\label{eqn:coef2}
\end{equation}
Therefore we obtain $C_{(i)k}$ by solving 
the linear equations (\ref{eqn:coef2}). 
Moreover, by using (\ref{eqn:coef1}), the noncommutative characteristic 
polynomial for the $i$-th row is written as 
\begin{eqnarray*}
\Phi_i(\lambda)&=&
\left| 
\begin{array}{cccc}
\cline{4-4}
\sum_{k=1}^n C_{(i)k}a_{i1}^{(n-k)} 
& \cdots 
& \sum_{k=1}^n C_{(i)k}a_{in}^{(n-k)} & 
\multicolumn{1}{|c|}{\lambda^n} \\ 
\cline{4-4}
a_{i1}^{(n-1)} & \cdots & a_{in}^{(n-1)} & \lambda^{n-1}\\
\vdots & & \vdots & \vdots \\
a_{i1}^{(0)} & \cdots & a_{in}^{(0)} & 1 \\
\end{array}
\right| \\
&& \\
&=&
\left| 
\begin{array}{cccc}
\cline{4-4}
0 & \cdots & 0 & 
\multicolumn{1}{|c|}{\lambda^n-\sum_{k=1}^n C_{(i)k}\lambda^{n-k}} \\ 
\cline{4-4}
a_{i1}^{(n-1)} & \cdots & a_{in}^{(n-1)} & \lambda^{n-1}\\
\vdots & & \vdots & \vdots \\
a_{i1}^{(0)} & \cdots & a_{in}^{(0)} & 1 \\
\end{array}
\right| \\
&=&\lambda^n-\sum_{k=1}^n C_{(i)k}\lambda^{n-k}.
\end{eqnarray*}
\end{proof}
By this proof, we obtain an important corollary. 
\begin{cor}
If an identity (\ref{eqn:NCCH}) holds, 
the noncommutative characteristic polynomials 
$\Phi_i(\lambda)$ defined by (\ref{eqn:NCCP}) 
are equal to $\lambda^n-\sum_{k=1}^n C_{(i)k}\lambda^{n-k}$. 
Especially, the (usual) Cayley-Hamilton theorem for $A$ 
(i.e. $C_{(i)k}=C_{k}$ for all $i$) holds, then 
$\Phi_i(\lambda) \ (i=1,\cdots,n)$ coincide with the usual 
characteristic polynomial $\Phi(\lambda)$ of $A$. 
\end{cor}
Moreover, as a contraposition, we have the following; 
\begin{cor}
For a given matrix $A$, if the noncommutative characteristic 
polynomials $\Phi_i(\lambda)$ are different for each $i$, then 
no commutative-Cayley-Hamilton-theorem type of identity with respect to $A$ 
exist. 
\end{cor}
\begin{exa}
Let $A$ be a matrix $A_q=(a_{ij})$ of the generators of the quantum group 
$GL_q(n)$, the noncommutative Cayley-Hamilton theorem 
(the quantum Cayley-Hamilton theorem) 
holds \cite{GKLLRT}. For example, $n=2$, by using relations 
$$ a_{11}a_{22}-a_{22}a_{11}=(q^{-1}-q)a_{12}a_{21}, \ \ 
a_{12}a_{22}=q^{-1}a_{22}a_{12}, \ \ 
a_{11}a_{21}=q^{-1}a_{21}a_{11}, \ \ 
a_{12}a_{21}=a_{21}a_{12}, $$
we have 
$$A_q^2-(q^{1/2}a_{11}+q^{-1/2}a_{22})
\begin{pmatrix}
q^{-1/2} & 0 \\
0 & q^{1/2}
\end{pmatrix}A_q
+(a_{11}a_{22}-q^{-1}a_{12}a_{21})
\begin{pmatrix}
q^{-1} & 0 \\
0 & q
\end{pmatrix}=O.$$
However, the noncommutative characteristic 
polynomials for each row do not coincide each other;
\begin{eqnarray*}
\Phi_1(\lambda)&=&\lambda^2-(a_{11}+q^{-1}a_{22})\lambda
+q^{-1}a_{11}a_{22}-q^{-2}a_{12}a_{21}, \\
\Phi_2(\lambda)&=&\lambda^2-(q \ a_{11}+a_{22})\ \lambda \ 
+q \ a_{11}a_{22}-a_{12}a_{21}.
\end{eqnarray*} 
Therefore, there is no identity for $A$ of 
commutative-Cayley-Hamilton-theorem type. 
\end{exa}

\section{Noncommutative Spectral Decomposition}

In this section, we develop a noncommutative analogue of 
the spectral decomposition with the quasideterminant. 
First, we review the Vandermonde quasideterminant and 
define a noncommutative analogue of the Lagrange interpolating polynomial. 
Next, we present the main theorem and our method of a noncommutative 
spectral decomposition. We also give a proof of the theorem 
by using properties of the quasideterminant prepared in section \ref{sec:pre}.

\subsection{Vandermonde Quasideterminant}

First, for $x_1, x_2, \cdots, x_k \in R \ $, 
the Vandermonde quasideterminant (\cite{GKLLRT}, \cite{GGRW}) is defined by 
$$V(x_1, \cdots, x_k)
=\left| 
\begin{array}{ccc}
\cline{3-3}
x_1^{k-1} & \cdots & \multicolumn{1}{|c|}{x_k^{k-1}} \\ \cline{3-3}
& \cdots & \\
x_1 & \cdots & x_k \\ 
1 & \cdots & 1 \\
\end{array}
\right|.
$$

\begin{exa}
\begin{eqnarray*}
V(x_1,x_2,z)&=&
\left| 
\begin{array}{ccc}
\cline{3-3}
x_1^2 & x_2^2 & \multicolumn{1}{|c|}{z^2} \\ \cline{3-3}
x_1 & x_2 & z \\ 
1 & 1 & 1 \\
\end{array}
\right| \\
&=&z^2-\begin{pmatrix}
x_1^2 & x_2^2 \\
\end{pmatrix}
\begin{pmatrix}
x_1 & x_2 \\
1 &  1
\end{pmatrix}^{-1}
\begin{pmatrix}
z \\ 1
\end{pmatrix}\\
&=&
z^2-\begin{pmatrix}
x_1^2 & x_2^2 \\
\end{pmatrix}
\begin{pmatrix}
(x_1-x_2)^{-1} & (1-x_2^{-1}x_1)^{-1} \\
(x_2-x_1)^{-1} & (1-x_1^{-1}x_2)^{-1} \\
\end{pmatrix}
\begin{pmatrix}
z \\ 1
\end{pmatrix}\\
&=&z^2+
\begin{pmatrix}
-x_1-(x_2-x_1)x_2(x_2-x_1)^{-1} & (x_2-x_1)x_2(x_2-x_1)^{-1}x_1 \\
\end{pmatrix}
\begin{pmatrix}
z \\ 1
\end{pmatrix}\\
&=&z^2+
\begin{pmatrix}
-(y_1+y_2) & y_2 y_1
\end{pmatrix}
\begin{pmatrix}
z \\ 1
\end{pmatrix}\\
&=&z^2-(y_1+y_2)z+y_2 y_1
\end{eqnarray*}
where we put $ y_1=x_1, \quad y_2=(x_2-x_1)x_2(x_2-x_1)^{-1}$. 
This is the noncommutative version of the relationship between solutions 
and coefficients for a (left) algebraic equation of degree $2$ \cite{GGRW}. 
\end{exa}

\begin{rem}
If $z=A=(a_{ij}) \in M(2,R)$ and $x_j=
\begin{pmatrix}
x_{(1)j} & \\
& x_{(2)j} 
\end{pmatrix} \ (j=1,2)$, $y_j \ (j=1,2)$ are also diagonal matrices. 
Moreover, comparing the equation $V(x_1,x_2,A)=A^2-(y_1+y_2)A+y_2 y_1=O$ 
with the noncommutative Cayley-Hamilton's theorem 
$$ A^2-
\begin{pmatrix}
C_{(1)1} & \\
& C_{(2)1} 
\end{pmatrix}A
-\begin{pmatrix}
C_{(1)2} & \\
& C_{(2)2} 
\end{pmatrix}=O, $$
if $y_1+y_2=\begin{pmatrix}
C_{(1)1} & \\
& C_{(2)1} 
\end{pmatrix}$
and $y_2 y_1=-\begin{pmatrix}
C_{(1)2} & \\
& C_{(2)2} 
\end{pmatrix}$, by the relationship between solutions and coefficients 
again, $x_{(i)1}, x_{(i)2}$ are the solutions of 
the noncommutative characteristic equation of $A$ for the $i$-th row. 

For a given $z=A=(a_{ij}) \in M(n,R)$ and the equation 
$V(x_1,\cdots,x_n,A)=O$, diagonal components of 
diagonal matrices $x_j$ are the solutions of the noncommutative 
characteristic equations of $A$ in the same way.
\label{rem:sol}
\end{rem}

\subsection{Noncommutative Lagrange Interpolating Polynomial}
For $x_1, \cdots, x_n \in R$, suppose that the inverse of 
the Vandermonde matrix 
$\displaystyle 
\begin{pmatrix}
x_1^{n-1} & \cdots & x_n^{n-1} \\
& \cdots & \\
1 & \cdots & 1
\end{pmatrix}^{-1}$ 
exist. Then we define polynomials 
$f_i(z) \ (i=1,\cdots,n)$ with respective to $z \in R$ as follows; 
\begin{defn}
$$\begin{pmatrix}
f_1(z) \\
\vdots \\
f_n(z)
\end{pmatrix}
=\begin{pmatrix}
x_1^{n-1} & \cdots & x_n^{n-1} \\
& \cdots & \\
1 & \cdots & 1
\end{pmatrix}^{-1}
\begin{pmatrix}
z^{n-1} \\
\vdots \\
1
\end{pmatrix}
$$
We call them noncommutative Lagrange interpolating polynomials. 
\end{defn}

\begin{exa}
For $n=2$, 
\begin{eqnarray*}
f_1(z)&=&\left| 
\begin{array}{cc}
\cline{1-1}
\multicolumn{1}{|c|}{x_1} & x_2 \\ \cline{1-1}
1 & 1 
\end{array}
\right|^{-1}z
+\left| 
\begin{array}{cc}
x_1 & x_2 \\ \cline{1-1}
\multicolumn{1}{|c|}{1} & 1 \\
\cline{1-1}
\end{array}
\right|^{-1} 1 \\
&=&(x_1-x_2)^{-1}z+(1-x_2^{-1}x_1)^{-1} \\
&=&(x_1-x_2)^{-1}z+(x_2-x_1)^{-1}x_2 \\
&=&(x_1-x_2)^{-1}(z-x_2),
\end{eqnarray*}
\begin{eqnarray*}
f_2(z)&=&\left| 
\begin{array}{cc}
\cline{2-2}
x_1 & \multicolumn{1}{|c|}{x_2}  \\ \cline{2-2}
1 & 1 
\end{array}
\right|^{-1}z
+\left| 
\begin{array}{cc}
x_1 & x_2 \\ \cline{2-2}
1 & \multicolumn{1}{|c|}{1} \\
\cline{2-2}
\end{array}
\right|^{-1} 1 \\
&=&(x_2-x_1)^{-1}z+(1-x_1^{-1}x_2)^{-1} \\
&=&(x_2-x_1)^{-1}(z-x_1).
\end{eqnarray*}
\end{exa}

By the definition above, we obtain the following theorem. 

\begin{thm}\ For $x_1, \cdots, x_n, \ z \in R$, we have \\
$(1) \quad x_1^j f_1(z)+x_2^j f_2(z)+\cdots+x_n^j f_n(z)=z^j 
\quad (j=0,1,\cdots,n-1)$, \\
$(2) \quad f_i(x_j)=\delta_{ij}.$
\end{thm}

\subsection{Our Method of Noncommutative Spectral Decomposition}

\begin{thm}[Main theorem]\label{thm:main}
If given $z, \ x_1, \cdots, x_n \in R$ satisfy the equation 
$V(x_1, \cdots, x_n, z)=0$, then we have the following identities 
\begin{equation}
V_m \equiv 
\left| 
\begin{array}{cccc}
\cline{4-4}
x_1^{m} & \cdots & x_n^{m} & \multicolumn{1}{|c|}{z^{m}} \\ \cline{4-4}
x_1^{n-1} & \cdots & x_n^{n-1} & z^{n-1}\\
& \cdots & & \\
1 & \cdots & 1 & 1 \\
\end{array}
\right|=0 \quad (m=0,\cdots,n,n+1,\cdots).
\label{eqn:zm}
\end{equation}
\end{thm}

Rewriting (\ref{eqn:zm}), by the definition of noncommutative 
Lagrange interpolating polynomials 
\begin{eqnarray*}
z^m
&=&
\begin{pmatrix}
x_1^{m} & \cdots & x_n^{m}
\end{pmatrix}
\begin{pmatrix}
x_1^{n-1} & \cdots & x_n^{n-1} \\
& \cdots & & \\
1 & \cdots & 1 
\end{pmatrix}^{-1}
\begin{pmatrix}
z^{n-1} \\
\vdots \\
1 
\end{pmatrix}
\\
&=&x_1^{m}f_1(z)+\cdots +x_n^{m}f_n(z), 
\end{eqnarray*}
then we have the noncommutative spectral decomposition of $z$
$$ z^m=x_1^{m}f_1(z)+\cdots +x_n^{m}f_n(z) 
\quad (m=0,1,\cdots).$$

In particular，if $z$ is a matrix $\ A=(a_{ij}) \in M(n,R)$, putting 
$x_1, \cdots , x_n$ as unknown diagonal matrices and solve the equation 
$$V(x_1,\cdots,x_n,A)=O. $$
By the remark \ref{rem:sol}, 
this equation is nothing but the noncommutative Cayley-Hamilton's 
theorem and the diagonal components of 
diagonal matrices $x_j$ are the solutions of the noncommutative 
characteristic equations of $A$. 
Therefore, by using the solutions of them, 
we obtain the noncommutative spectral decomposition of $A$
$$ A^m=x_1^{m}f_1(A)+\cdots +x_n^{m}f_n(A) 
\quad (m=0,1,\cdots). $$

\subsection{A Proof of Main Theorem \ref{thm:main}}

\begin{proof}
In case of $m=0,1,\cdots,n-1$, the identity (\ref{eqn:zm}) is trivial. 
If $m=n$, (\ref{eqn:zm}) is nothing but $V(x_1, \cdots, x_n, z)=0$. 
In the following, we suppose $m=n+1,n+2,\cdots$. 

Consider a matrix $A$ and the submatrix $A_0$ defined by 
$$A=
\left(
\begin{array}{ccc|cc}
x_1^{m} & \cdots & x_n^{m} & 0 & z^{m} \\
x_1^{n} & \cdots & x_n^{n} & 0 & z^{n} \\ \cline{1-5}
x_1^{n-1} & \cdots & x_n^{n-1} & 0 & z^{n-1}\\
& \cdots & & \vdots & \vdots\\
1 & \cdots & 1 & 1 & 1 \\
\end{array}
\right), 
\qquad 
A_0=
\begin{pmatrix}
x_1^{n-1} & \cdots & x_n^{n-1} \\
& \cdots & \\
1 & \cdots & 1 \\
\end{pmatrix}.$$
For $p=1,2, \ q=n+1,n+2$, we put a matrix $C=(c_{pq})$ which entries are 
quasideterminants with $A_0$ as a pivot like (\ref{eqn:pivot}) 
(note that quasideterminants are unchanged under permutations of 
rows or columns) and we remark 
$$\hspace{-31mm}
c_{1,n+2}=
\left|
\begin{array}{cccc}
\cline{4-4}
x_1^{m} & \cdots & x_n^{m} & \multicolumn{1}{|c|}{z^{m}} \\ \cline{4-4}
&  &  & z^{n-1}\\ 
& A_0 & & \vdots\\
&  &  & 1 \\
\end{array}
\right|=V_m, $$
$$
c_{2,n+2}=
\left|
\begin{array}{cccc}
\cline{4-4}
x_1^{n} & \cdots & x_n^{n} & \multicolumn{1}{|c|}{z^{n}} \\ \cline{4-4}
&  &  & z^{n-1}\\ 
& A_0 & & \vdots\\
&  &  & 1 \\
\end{array}
\right|=V(x_1, \cdots, x_n, z)=V_n.
$$
Then by the Sylvester's identity
(theorem \ref{thm:Syl}), we have 
\begin{eqnarray*}
|A|_{1,n+2}
&=&|C|_{1,n+2}
=\left|
\begin{array}{cc}
\cline{2-2}
c_{1,n+1} & \multicolumn{1}{|c|}{c_{1,n+2}} \\ \cline{2-2}
c_{2,n+1} & c_{2,n+2} \\
\end{array}
\right| \\
&=&c_{1,n+2}-c_{1,n+1}c_{2,n+1}^{-1}c_{2,n+2} \\
&=&V_m-c_{1,n+1}c_{2,n+1}^{-1}V_n. 
\end{eqnarray*}
On the other hand, since 
\begin{eqnarray*}
|A|_{1,n+2}
&=&\left|
\begin{array}{ccc|cc}
\cline{5-5}
x_1^{m} & \cdots & x_n^{m} & 0 & \multicolumn{1}{|c|}{z^{m}} \\ \cline{5-5}
x_1^{n} & \cdots & x_n^{n} & 0 & z^{n} \\ \cline{1-5}
x_1^{n-1} & \cdots & x_n^{n-1} & 0 & z^{n-1}\\
& \cdots & & \vdots & \vdots\\
1 & \cdots & 1 & 1 & 1 \\
\end{array}
\right| 
=\left|
\begin{array}{cccc|c}
\cline{4-4}
x_1^{m} & \cdots & x_n^{m} & \multicolumn{1}{|c|}{z^{m}} & 0 \\ \cline{4-4}
x_1^{n} & \cdots & x_n^{n} & z^{n} & 0 \\ 
x_1^{n-1} & \cdots & x_n^{n-1} & z^{n-1} & 0 \\
& \cdots & & \vdots & \vdots\\ \cline{1-5}
1 & \cdots & 1 & 1 & 1 \\
\end{array}
\right|\\
&& \\
&=&\left|
\begin{array}{cccc}
\cline{4-4}
x_1^{m} & \cdots & x_n^{m} & \multicolumn{1}{|c|}{z^{m}} \\ \cline{4-4}
x_1^{n} & \cdots & x_n^{n} & z^{n} \\ 
x_1^{n-1} & \cdots & x_n^{n-1} & z^{n-1} \\
& \cdots & & \vdots \\ 
x_1 & \cdots & x_n & z \\
\end{array}
\right| 
\quad
\begin{pmatrix}
\mbox{by the Sylvester's identity with}\\
\mbox{$(n+2,n+2)$-entry $1$ as a pivot}
\end{pmatrix} \\
&& \\
&=&\left|
\begin{array}{cccc}
\cline{4-4}
x_1^{m-1} & \cdots & x_n^{m-1} & \multicolumn{1}{|c|}{z^{m-1}} \\ 
\cline{4-4}
x_1^{n-1} & \cdots & x_n^{n-1} & z^{n-1} \\ 
x_1^{n-2} & \cdots & x_n^{n-2} & z^{n-2} \\
& \cdots & & \vdots \\ 
1 & \cdots & 1 & 1 \\
\end{array}
\right| z 
\quad 
\begin{pmatrix}
\mbox{by the property (\ref{eqn:mu})}
\end{pmatrix} \\
&=& V_{m-1}z, 
\end{eqnarray*}
we obtain the important identities
\begin{equation}
V_{m-1}z=V_m-c_{1,n+1}c_{2,n+1}^{-1}V_n.
\label{eqn:main-identity}
\end{equation}
Therefore $V_n=0$ implies $V_m=0 \ (m=n+1,n+2,\cdots)$ by 
the mathematical induction. 
\end{proof}

\begin{rem}
In particular, for $n=2$, 
$$V_m=
z^m-
\begin{pmatrix}
x_1^{m} & x_2^{m}
\end{pmatrix}
\begin{pmatrix}
x_1 & x_2 \\
1 & 1 
\end{pmatrix}^{-1}
\begin{pmatrix}
z \\
1 
\end{pmatrix}
=z^m-(x_1^{m}f_1(z)+x_2^{m}f_2(z)).
$$
Then the identity (\ref{eqn:main-identity}) is 
\begin{eqnarray}
&&\{ z^{m-1}-(x_1^{m-1}f_1(z)+x_2^{m-1}f_2(z))\} z \nonumber\\
&=& z^m-(x_1^{m}f_1(z)+x_2^{m}f_2(z)) \nonumber\\
&& \ -(x_2^{m-1}-x_1^{m-1})(x_2-x_1)^{-1}
\{ z^2-(x_1^{2}f_1(z)+x_2^{2}f_2(z)) \} 
\qquad (m=2,3,\cdots).
\end{eqnarray}
\end{rem}


\section{Examples of Noncommutative Spectral Decomposition and 
the Exponential Matrices}

In this section, we apply our method to a quaternionic matrix and 
a matrix with entries being harmonic oscillators. As a result, 
we obtain the noncommutative spectral decomposition and 
the exponential matrices of them explicitly. 

\subsection{Quaternionic Matrices}

As a quaternionic matrix, we consider an element $A$ of 
Lie algebra $sp(2)$; 
$$A=
\begin{pmatrix}
i & j \\
j & -i 
\end{pmatrix}
$$
We apply our method to $A$ and calculate the spectral decomposition 
and the exponential matrix $\exp tA$ explicitly. First, from 
$$A^2=\begin{pmatrix}
-2 & 2k \\
-2k & -2
\end{pmatrix},
$$
noncommutative characteristic equations for each row are 
$$
\Phi_1(\lambda)
=
\left| 
\begin{array}{cccc}
\cline{3-3}
-2 & 2k & \multicolumn{1}{|c|}{\lambda^2} \\ \cline{3-3}
i & j & \lambda \\
1 & 0 & 1 \\
\end{array}
\right|
=\lambda^2-2i\lambda=0, 
$$
then $\lambda=0, \ 2i$, and 
$$
\Phi_2(\lambda)
=
\left| 
\begin{array}{cccc}
\cline{3-3}
-2k & -2 & \multicolumn{1}{|c|}{\lambda^2} \\ \cline{3-3}
j & -i & \lambda \\
0 & 1 & 1 \\
\end{array}
\right|
=\lambda^2+2i\lambda=0, 
$$
then $\lambda=0, \ -2i$.

Next, in the noncommutative 
Lagrange interpolating polynomials $$
f_1(z)=(x_1-x_2)^{-1}(z-x_2), \quad
f_2(z)=(x_2-x_1)^{-1}(z-x_1), 
$$
we put 
$z=A, \quad 
x_1=\begin{pmatrix}
2i & \\
 & -2i
\end{pmatrix}, \quad 
x_2=\begin{pmatrix}
0 & \\
 & 0
\end{pmatrix}
$, then 
\begin{eqnarray*}
P_1&=&f_1(A)=
\begin{pmatrix}
2i & \\
 & -2i
\end{pmatrix}^{-1}A
=\frac{1}{2}\begin{pmatrix}
1 & -k \\
k & 1
\end{pmatrix}, \\
P_2&=&f_2(A)=
\left\{ 
-\begin{pmatrix}
2i & \\
 & -2i
\end{pmatrix}
\right\}^{-1}
\left\{ A-
\begin{pmatrix}
2i & \\
 & -2i
\end{pmatrix}
\right\}
=\frac{1}{2}\begin{pmatrix}
1 & k \\
-k & 1
\end{pmatrix}.
\end{eqnarray*}
We can check 
$\ P_i^2=P_i, \ P_iP_j=0 \ (i \ne j)$ easily and we obtain 
\begin{eqnarray*}
\exp tA &=& (\exp tx_1) P_1+(\exp tx_2) P_2 \\
&=& \begin{pmatrix}
e^{2it} & \\
 & e^{-2it}
\end{pmatrix}
\frac{1}{2}\begin{pmatrix}
1 & -k \\
k & 1
\end{pmatrix}
+\begin{pmatrix}
1 &  \\
 & 1
\end{pmatrix}
\frac{1}{2}\begin{pmatrix}
1 & k \\
-k & 1
\end{pmatrix}\\
&=&\frac{1}{2}\begin{pmatrix}
e^{2it}+1 & -e^{2it}k+k \\
e^{-2it}k-k & e^{-2it}+1
\end{pmatrix}\\
&=&\begin{pmatrix}
e^{it}\cos t & e^{it}\sin t \ j \\
e^{-it}\sin t \ j & e^{-it}\cos t
\end{pmatrix} \ \in Sp(2). 
\end{eqnarray*}

\begin{rem}
If we put 
$$x_1=\begin{pmatrix}
2i & \\
 & 0
\end{pmatrix}, \quad 
x_2=\begin{pmatrix}
0 & \\
 & -2i
\end{pmatrix},
\ $$
we have $[A, x_1] \ne O, \ [A, x_2] \ne O$ and 
$f_1(A), \ f_2(A)$ are {\textbf{not}} projection matrices. 
Nevertheless, by the theorem \ref{thm:main}, we have 
$A^m=x_1^m f_1(A)+x_2^m f_2(A)$ and we can calculate 
$\exp tA$ explicitly！ This result is derived from the fact 
that the theorem \ref{thm:main} 
is not depend on the ordering of solutions for 
noncommutative characteristic equations for each row. 
\end{rem}

\subsection{A Matrix with Entries being Harmonic Oscillators}

Let $a,a^{\dagger}$ be the generator of the harmonic oscillator．
The relation is $ \ [a,a^{\dagger}]=1.$
We also denote $N$ as the number operator $N=a^{\dagger}a$. 

We consider a matrix 
$\displaystyle 
A=\sqrt{2}
\begin{pmatrix}
0 & a & 0 \\
\adag & 0 & a \\
0 & \adag & 0 \\
\end{pmatrix}$. 
This matrix is related to a Hamiltonian of a model 
in quantum optics \cite{FHKSW}. So, it is important to calculate the 
exponential of $A$ as the time-evolution operator of the Hamiltonian. 

From 
$$
A^2=
2\begin{pmatrix}
N+1 & 0 & a^2 \\
0 & 2N+1 & 0 \\
(\adag)^2 & 0 & N \\
\end{pmatrix}, \quad 
A^3=2\sqrt{2}
\begin{pmatrix}
0 & (2N+3)a & 0 \\
(2N+1)\adag & 0 & (2N+1)a \\
0 & (2N-1)\adag & 0 \\
\end{pmatrix}, 
$$
the noncommutative characteristic equations for each row are 
$$
\Phi_1(\lambda)
=
\left| 
\begin{array}{cccc}
\cline{4-4}
0 & 2\sqrt{2}(2N+3)a & 0 & \multicolumn{1}{|c|}{\lambda^3} \\ \cline{4-4}
2(N+1) & 0 & 2a^2 & \lambda^2 \\
0 & \sqrt{2}a & 0 & \lambda \\ 
1 & 0 & 0 & 1 \\
\end{array}
\right|
=\lambda^3-2(2N+3)\lambda=0, 
$$
then $\lambda=\pm \sqrt{2(2N+3)}, \ 0$, and 
$$
\Phi_3(\lambda)
=
\left| 
\begin{array}{cccc}
\cline{4-4}
0 & 2\sqrt{2}(2N-1)\adag & 0 & \multicolumn{1}{|c|}{\lambda^3} \\ \cline{4-4}
2(\adag)^2 & 0 & 2N & \lambda^2 \\
0 & \sqrt{2}\adag & 0 & \lambda \\ 
0 & 0 & 1 & 1 \\
\end{array}
\right|
=\lambda^3-2(2N-1)\lambda=0,
$$
then $\lambda=\pm \sqrt{2(2N-1)}, \ 0$.

\begin{rem}
For the second row, the quasideterminant 
$$
\Phi_2(\lambda)
=
\left| 
\begin{array}{cccc}
\cline{4-4}
2\sqrt{2}(2N+1)\adag & 0 & 2\sqrt{2}(2N+1)a & 
 \multicolumn{1}{|c|}{\lambda^3} \\ \cline{4-4}
0 & 2(2N+1) & 0 & \lambda^2 \\
\sqrt{2}\adag & 0 & \sqrt{2}a & \lambda \\ 
0 & 1 & 0 & 1 \\
\end{array}
\right|
$$
is not defined because ``rank" of the matrix 
$
\begin{pmatrix}
0 & 2(2N+1) & 0 \\
\sqrt{2}\adag & 0 & \sqrt{2}a \\ 
0 & 1 & 0
\end{pmatrix}$
is $2$ (on ``rank of $A \in M(n,R)$", see \cite{GGRW}). Then, we put 
$$A^3+UA^2+VA+W=0, \quad U,V,W \mbox{ are diagonal matrices}$$
and simplify them, then the second row is 
$$\lambda^3+u \lambda^2-2(2N+1)\lambda-2u(2N+1)=0 \quad 
(\mbox{for arbitrary } u). $$
Therefore, if we put $u=0$, we obtain $\lambda=\pm \sqrt{2(2N+1)}, \ 0$. 
\end{rem}

Next, we calculate the noncommutative 
Lagrange interpolating polynomials 
\begin{eqnarray*}
f_1(z)&=&\left| 
\begin{array}{ccc}
\cline{1-1}
\multicolumn{1}{|c|}{x_1^2} & x_2^2 & x_3^2 \\ \cline{1-1}
x_1 & x_2 & x_3 \\
1 & 1 & 1
\end{array}
\right|^{-1}z^2
+\left| 
\begin{array}{ccc}
x_1^2 & x_2^2 & x_3^2 \\
\cline{1-1}
\multicolumn{1}{|c|}{x_1} & x_2 & x_3 \\ \cline{1-1}
1 & 1 & 1
\end{array}
\right|^{-1}z
+\left| 
\begin{array}{ccc}
x_1^2 & x_2^2 & x_3^2 \\
x_1 & x_2 & x_3 \\ \cline{1-1}
\multicolumn{1}{|c|}{1} & 1 & 1 \\
\cline{1-1}
\end{array}
\right|^{-1} 1 \\
&& \\
&=& \ 
(x_1^2-x_2^2(x_2-x_3)^{-1}x_1-x_3^2(x_3-x_2)^{-1}x_1 \\
&& \hspace{2cm}
-x_2^2(x_3-x_2)^{-1}x_3-x_3^2(x_2-x_3)^{-1}x_2)^{-1}z^2 \\
&&+(x_1-x_2(x_2^2-x_3^2)^{-1}x_1^2-x_3(x_3^2-x_2^2)^{-1}x_1^2 \\
&& \hspace{2cm}
-x_2(x_3^2-x_2^2)^{-1}x_3^2-x_3(x_2^2-x_3^2)^{-1}x_2^2)^{-1}z \\
&&+(1-(x_2^2-x_3x_2)^{-1}x_1^2-(x_3^2-x_2x_3)^{-1}x_1^2 \\
&& \hspace{2cm}
-(x_2-x_3^{-1}x_2^2)^{-1}x_1-(x_3-x_2^{-1}x_3^2)^{-1}x_1)^{-1}. 
\end{eqnarray*}
In particular, in the case of $x_1=x, \ x_2=0, \ x_3=-x $,
\begin{eqnarray*}
f_1(z)&=&(x^2-(-x)^2(-x)^{-1}x)^{-1}z^2+(x-(-x)(-x)^{-2}x^2)^{-1}z+0\\
&=&(2x^2)^{-1}z^2+(2x)^{-1}z, 
\end{eqnarray*}
where the last term of $f_1(z)$ is calculated by using the homological 
relation as follows; 
\begin{eqnarray*}
\left| 
\begin{array}{ccc}
x_1^2 & x_2^2 & x_3^2 \\
x_1 & x_2 & x_3 \\ \cline{1-1}
\multicolumn{1}{|c|}{1} & 1 & 1 \\
\cline{1-1}
\end{array}
\right|^{-1}
&=&-\left| 
\begin{array}{ccc}
x_1^2 & x_2^2 & x_3^2 \\
\cline{1-1}
\multicolumn{1}{|c|}{x_1} & x_2 & x_3 \\ \cline{1-1}
1 & 1 & 1
\end{array}
\right|^{-1}
\left| 
\begin{array}{cc}
x_2^2 & x_3^2 \\ \cline{1-1}
\multicolumn{1}{|c|}{x_2} & x_3 \\ \cline{1-1}
\end{array}
\right|
\left| 
\begin{array}{cc}
x_2^2 & x_3^2 \\ \cline{1-1}
\multicolumn{1}{|c|}{1} & 1 \\ \cline{1-1}
\end{array}
\right|^{-1}\\
&=&
-\left| 
\begin{array}{ccc}
x_1^2 & x_2^2 & x_3^2 \\
\cline{1-1}
\multicolumn{1}{|c|}{x_1} & x_2 & x_3 \\ \cline{1-1}
1 & 1 & 1
\end{array}
\right|^{-1}
(x_2-x_3^{-1}x_2^2)(1-x_3^{-2}x_2^2)^{-1}\\
&\to & -(2x)^{-1} \cdot 0 \cdot 1=0 
\quad (x_1 \to x, \ x_2 \to 0, \ x_3 \to -x).
\end{eqnarray*}
Then we put 
$$z=A, \quad 
x=
\begin{normalsize}
\begin{pmatrix}
\lambda(N) & & \\
& \lambda(N-1) & \\
& & \lambda(N-2) \\
\end{pmatrix}
\end{normalsize}, \ 
\lambda(N)=\sqrt{2(2N+3)}, $$
we have 
\begin{eqnarray*}
P_1&=&f_1(A)\\
&=&(2x^2)^{-1}A^2+(2x)^{-1}A \\
&=&\begin{normalsize}
\begin{pmatrix}
(2\lambda(N))^{-2} & & \\
& (2\lambda(N-1))^{-2} & \\
& & (2\lambda(N-2))^{-2} \\
\end{pmatrix}\cdot 
2\begin{pmatrix}
N+1 & 0 & a^2 \\
0 & 2N+1 & 0 \\
(\adag)^2 & 0 & N \\
\end{pmatrix}
\end{normalsize}\\
&& \\
&&\hspace{1cm}+\begin{normalsize}
\begin{pmatrix}
(2\lambda(N))^{-1} & & \\
& (2\lambda(N-1))^{-1} & \\
& & (2\lambda(N-2))^{-1} \\
\end{pmatrix}\cdot 
\sqrt{2}\begin{pmatrix}
0 & a & 0 \\
\adag & 0 & a \\
0 & \adag & 0 \\
\end{pmatrix}
\end{normalsize}\\
&& \\
&& \\
&=&\begin{large}
\begin{pmatrix}
\frac{1}{2(2N+3)} & & \\
& \frac{1}{2(2N+1)} & \\
& & \frac{1}{2(2N-1)} \\
\end{pmatrix}
\begin{pmatrix}
N+1 & 0 & a^2 \\
0 & 2N+1 & 0 \\
(\adag)^2 & 0 & N \\
\end{pmatrix}
\end{large}\\
&& \\
&&\hspace{1cm}+\begin{large}
\begin{pmatrix}
\frac{1}{2\sqrt{2N+3}} & & \\
& \frac{1}{2\sqrt{2N+1}} & \\
& & \frac{1}{2\sqrt{2N-1}} \\
\end{pmatrix}
\begin{pmatrix}
0 & a & 0 \\
\adag & 0 & a \\
0 & \adag & 0 \\
\end{pmatrix}
\end{large}\\
&& \\
&& \\
&=&
\begin{pmatrix}
\frac{N+1}{2(2N+3)} & \frac{1}{2\sqrt{2N+3}}a & 
\frac{1}{2(2N+3)}a^2 \\
\frac{1}{2\sqrt{2N+1}}\adag & \frac{1}{2} & \frac{1}{2\sqrt{2N+1}}a \\
\frac{1}{2(2N-1)}(\adag)^2 & \frac{1}{2\sqrt{2N-1}}\adag & 
\frac{N}{2(2N-1)}
\end{pmatrix}.
\end{eqnarray*}

In the same manner, if we put 
$z=A, \ x_1=x, \ x_2=0, \ x_3=-x$ in $f_2(z), \ f_3(z)$, then we have 
\begin{eqnarray*}
P_2&=&f_2(A)\\
&=&(-x^2)^{-1}A^2+I_3 \\
&=&\begin{large}
-\begin{pmatrix}
\frac{1}{2(2N+3)} & & \\
& \frac{1}{2(2N+1)} & \\
& & \frac{1}{2(2N-1)} \\
\end{pmatrix}\cdot 
2\begin{pmatrix}
N+1 & 0 & a^2 \\
0 & 2N+1 & 0 \\
(\adag)^2 & 0 & N \\
\end{pmatrix}
+\begin{pmatrix}
1 & & \\
& 1 & \\
& & 1 \\
\end{pmatrix}
\end{large}\\
&& \\
&=&
\begin{pmatrix}
\frac{N+2}{2N+3} & 0 & -\frac{1}{2N+3}a^2 \\
0 & 0 & 0 \\
-\frac{1}{2N-1}(\adag)^2 & 0 & \frac{N-1}{2N-1}
\end{pmatrix},
\end{eqnarray*}
\vspace{2mm}
\begin{eqnarray*}
P_3=f_3(A)&=&(2x^2)^{-1}A^2-(2x)^{-1}A \\
&=&\begin{large}
\begin{pmatrix}
\frac{1}{2(2N+3)} & & \\
& \frac{1}{2(2N+1)} & \\
& & \frac{1}{2(2N-1)} \\
\end{pmatrix}
\begin{pmatrix}
N+1 & 0 & a^2 \\
0 & 2N+1 & 0 \\
(\adag)^2 & 0 & N \\
\end{pmatrix}
\end{large}\\
&& \\
&&\hspace{1cm}-\begin{large}
\begin{pmatrix}
\frac{1}{2\sqrt{2N+3}} & & \\
& \frac{1}{2\sqrt{2N+1}} & \\
& & \frac{1}{2\sqrt{2N-1}} \\
\end{pmatrix}
\begin{pmatrix}
0 & a & 0 \\
\adag & 0 & a \\
0 & \adag & 0 \\
\end{pmatrix}
\end{large}\\
&& \\
&=&
\begin{pmatrix}
\frac{N+1}{2(2N+3)} & -\frac{1}{2\sqrt{2N+3}}a & 
\frac{1}{2(2N+3)}a^2 \\
-\frac{1}{2\sqrt{2N+1}}\adag & \frac{1}{2} & -\frac{1}{2\sqrt{2N+1}}a \\
\frac{1}{2(2N-1)}(\adag)^2 & -\frac{1}{2\sqrt{2N-1}}\adag & 
\frac{N}{2(2N-1)}
\end{pmatrix}.
\end{eqnarray*}
We can check 
$\ P_i^2=P_i, \ P_iP_j=0 \ (i \ne j)$ and for a constant $g$, 
\begin{eqnarray*}
&& \exp(-itgA) \nonumber\\
&=& \exp(-itgx) P_1+(\exp 0) P_2+\exp(itgx) P_3 \nonumber\\
&=& \hspace{-3mm}
\begin{normalsize}
\begin{pmatrix}
\frac{N+2+(N+1)\cos(tg\lambda(N))}{2N+3} & 
-i\frac{1}{\sqrt{2N+3}}\sin(tg\lambda(N))a & 
\frac{1}{2N+3}(-1+\cos(tg\lambda(N)))a^2 \\
-i\frac{1}{\sqrt{2N+1}}\sin(tg\lambda(N-1))\adag & 
\cos(tg\lambda(N-1)) & 
-i\frac{1}{\sqrt{2N+1}}\sin(tg\lambda(N-1))a \\
\frac{1}{2N-1}(-1+\cos(tg\lambda(N-2)))(\adag)^2 & 
-i\frac{1}{\sqrt{2N-1}}\sin(tg\lambda(N-2))\adag & 
\frac{N-1+N\cos(tg\lambda(N-2))}{2N-1}
\end{pmatrix}.
\end{normalsize}
\end{eqnarray*}
\begin{rem}
For this $A$, by using ``the quantum diagonalization method" \cite{FHKSW}, 
$$\begin{normalsize}
\begin{pmatrix}
1 & &  \\
  & a\frac{1}{\sqrt{N}} &  \\
  & & a^2 \frac{1}{\sqrt{N(N-1)}}
\end{pmatrix}
A 
\begin{pmatrix}
1 & &  \\
  & \frac{1}{\sqrt{N}} \adag &  \\
  & & \frac{1}{\sqrt{N(N-1)}} (\adag)^2  
\end{pmatrix}
=\sqrt{2}
\begin{pmatrix}
0 & \sqrt{N+1} & 0 \\
\sqrt{N+1} & 0 & \sqrt{N+2} \\
0 & \sqrt{N+2} & 0 
\end{pmatrix}.
\end{normalsize}$$
Since the matrix on the right-hand side has only commutative entries, 
we calculate the characteristic equation as usual, then 
$$ \lambda^3-2(2N+3)\lambda=0, \quad 
\quad \lambda=0, \ \pm \sqrt{2(2N+3)}.$$
We remark that the result of $\exp(-itgA)$ in \cite{FHKSW} and with our 
noncommutative spectral decomposition coinside. 
\end{rem}

\section{Discussion}
In this paper, we developed a noncommutative version of the spectral 
decomposition with the quasideterminant 
and calculated some interesting examples. 
In particular, we defined a noncommutative analogue of 
the Lagrange interpolating polynomials and applied to the systematic method 
for constructing projection matrices with noncommutative entries. 

Our method is very powerful to calculate a function of a matrix with 
noncommutative entries and is expected to apply for the theory of 
noncommutative geometry, quantum physics, and so on. A study of 
other applications with our theory is in progress. 
 \section*{Acknowledgements}
The author is very grateful to Kazuyuki Fujii for helpful comments 
on an earlier draft on this paper and to Masashi Hamanaka and 
Hideyuki Ishi for helpful suggestion. 

\end{document}